\documentclass[11pt]{article}

\usepackage{graphicx}

\newtheorem{thm}{Theorem}[section]

\begin{document}

%

\title{{\LARGE\bf A Paradox for Admission Control of Multiclass Queueing
Network with Differentiated Service}}

\author{
Heng-Qing Ye\\
School of Business, National University of Singapore \\
1 Business Link,Singapore, bizyehq@nus.edu.sg }
\date{July 2005}

\maketitle

\begin{abstract}
In this paper, we present counter-intuitive examples for the
multiclass queueing network system. In the system, each station may
serve more than one job class with differentiated service priority,
and each job may require service sequentially by more than one
service station. In our examples, the network performance is
improved even when more workloads are admitted for service. \\

\end{abstract}

\section{Introduction}

The queueing network model is an important tool in studying the
service system, the manufacturing system and the communication
system. In many applications, the model is useful in identifying
bottleneck service resources of a system so that better decisions
can be made on designing and controlling the network. Such decisions
may include, for example, selecting the system service capacity
(e.g., the maximum service rates of work stations), adjusting system
workload (e.g., the job arrival rate and pattern), and routing jobs
to service stations if the jobs can be served by more than one
station.

In practice, it is commonly believed that the performance for a
queueing network system, say in terms of the average total number or
the average delay of jobs in the system, would be improved if the
service capacity (system workload, routing alternatives, resp.) is
increased (decreased, increased, resp.). Such an understanding is
sound when studying, e.g., the queuing system with single or
parallel service stations and the product-form queueing network; cf.
Chen and Yao (2001) and references there.
 However, one must be cautious in applying such an intuition to
complex queueing systems. In fact, from the study of the stability
condition of a three-station multiclass queueing network in Dumas
(1997), it can be noticed that the network (Dumas network) with
increased service capacity for certain work stations performs worse.
The paradox about the (distributed) routing in the queueing network
can also be found in Cohen and Kelly (1990), which is based on the
well known Braess paradox (Braess 1968).
 Complement to these paradoxes on the service capacity and
routing, we provide a paradoxical network examples on the admission
control. These counter-intuitive examples show that the network
performance could be degraded even when the system workload or the
arrival rate of jobs decreases.

We describe the multiclass queueing network model and present the
counter-intuitive results in the next section. In Section 3, we
introduce the fluid model approach developed in recent years and
then use this approach to prove our main results. We conclude in
Section 4.

\section{Counter-examples and Main Results} \label{sec-main}


The {\it multiclass queueing network} consists of $J$ stations
indexed by $j\in {\cal J} = \{ 1,\cdots,J\},$ and $K$ job classes
indexed by $k \in {\cal K} = \{1,\cdots,K\}.$ Assume that the
arrival process of class $k$ jobs (or customers) is a Poisson
process with arrival rate $\alpha_k~(\ge0)$, and the service time
for each class $k$ job is exponentially distributed with mean
service time $m_k~(>0)$. Denote
$\alpha=(\alpha_1,\cdots,\alpha_J)^T$ and $m=(m_1,\cdots,m_K)^T$. We
also assume that all the interarrival times and service times are
independent. A class $k$ job is served at station $\sigma(k)$
($\sigma(\cdot):~{\cal K}\rightarrow {\cal J}$), and after its
service completion, it may become a class $\ell$ job with
probability $p_{k\ell}$ and leave the network with probability
$1-\sum_{\ell=1}^Kp_{k\ell}$. Let $P=(p_{k\ell})$. Let $C=(c_{jk})$
be a $J\times K$ matrix whose $(j,k)$th component $c_{jk}=1$ if
$j=\sigma(k)$ and $=0$ otherwise. While each station may serve more
than one class of jobs, each job is served at one specific station
(determined by the many-to-one mapping $\sigma(\cdot)$). We study a
preemptive priority service discipline in this paper. Let $\pi$ be a
one-to-one mapping from ${\cal K}$ onto ${\cal K}$. For any given
$\ell$ and $k$, if $\pi(\ell)<\pi(k)$ and $\sigma(\ell)=\sigma(k)$,
then class $k$ job can not be served at station $\sigma(k)$ unless
there is no class $\ell$ job. In short, we say that class $\ell$ has
a higher priority than class $k$.  For convenience, the mapping
$\pi$ is often expressed as a permutation of ${\cal K}$, i.e., which
can be written as $\pi=(i_1,\cdots,i_K)$ if $\pi(k)=i_k$, $k\in{\cal
K}$. In addition, we only consider work-conserving (or non-idling)
service disciplines, which specify that a work station can not be
idle unless there is no job waiting for service in that station. For
convenience, we denote the queueing network described above as
(${\cal J}, {\cal K}, \alpha, m, C, P, \pi$).

We study {\it open} multiclass queueing network in this paper, or we
assume that the transition $P$ is transient, i.e.,
\begin{eqnarray}
&& I+P+P^2+\cdots \mbox{~~is convergent.} \label{eq-sn-80}
\end{eqnarray}
Let $\lambda=(I-P')^{-1}\alpha$, $\beta=M\lambda$ and
$\rho=C\beta=CM\lambda$. Call $\lambda$ a {\it nominal total
arrival rate\/} (vector), $\beta_k$ (the $k$th component of
$\beta$) a {\it traffic intensity for class\/} $k$, $k\in {\cal
K}$, and $\rho_j$ (the $j$th component of $\rho$) a {\it traffic
intensity for station\/} $j$, $j\in {\cal J}$. Usually, the vector
$\rho=(\rho_j)$ is simply called the {\it traffic intensity of the
queueing network}. Actually, $\lambda$ is the unique solution to
the following {\it traffic equation},
    $$\lambda = \alpha + P' \lambda,$$
which indicates that the nominal total arrival rate vector
$\lambda$ includes both external arrivals and internal
transitions.

The dynamics of the network can be described using a $K$-dimensional
queue length process $Q(t)=(Q_k(t), k\in {\cal K})$ ($t\ge 0$),
where $Q_k(t)$ indicates the number of class $k$ jobs in the network
at time $t$. The queue length process $Q(t)$ is a continuous time
Markov chain under the Poisson arrival and exponential service
assumptions. We say that the network (${\cal J}, {\cal K}, \alpha,
m, C, P, \pi$) is stable if the Markov chain $Q(t)$ is positive
recurrent. It is well know that the Markov chain $Q(t)$ is positive
recurrent {\it only if} the traffic intensity for each station is
less than one, i.e., $\rho_j < 1$ for all $j\in {\cal J}$, or in
short, $\rho < e$ where $e$ is a $J$-dimensional vector with all
components being ones. The performance index of interest in this
paper is the expected stationary total queue length $\bar Q$ defined
as
    $$\bar Q = \lim_{t\rightarrow \infty} E \left[ \sum_{k\in {\cal K}} Q_k(t) \right]. $$
The queue length $\bar Q(t)$ is a finite if and only if the queue
length process $Q$ is positive recurrent.

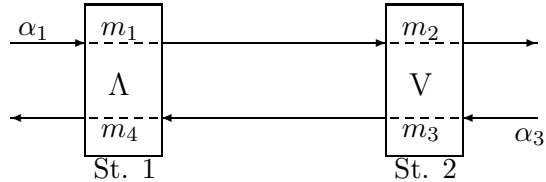
\begin{figure}[htb]
\setlength{\unitlength}{1cm}
\begin{center}
\begin{picture}(5,2)(1,1.5)
\multiput(1,1.5)(4,0){2}{\framebox(1,2)} 
\put(0,3){\vector(1,0){1}} 
\put(2,3){\vector(1,0){3}} 
\put(6,3){\vector(1,0){1}} 
   \multiput(1.05,3)(0.2,0){5}{\line(1,0){0.1}}
   \multiput(5.05,3)(0.2,0){5}{\line(1,0){0.1}}
\multiput(1,2)(6,0){2}{\vector(-1,0){1}} 
\put(5,2){\vector(-1,0){3}} 
   \multiput(1.05,2)(0.2,0){5}{\line(1,0){0.1}}
   \multiput(5.05,2)(0.2,0){5}{\line(1,0){0.1}}
\put(1.1,1.2){St. 1} \put(5.1,1.2){St. 2}
\put(0.1,3.1){$\alpha_1$}\put(6.7,1.7){$\alpha_3$}%
\put(1.2,3.1){$m_1$} \put(5.2,3.1){$m_2$} \put(1.2,1.75){$m_4$}
\put(5.2,1.75){$m_3$} \put(1.3,2.3){\large $\Lambda$}
\put(5.3,2.3){\large V}
\end{picture}
\end{center}
\caption{KRSS network} \label{fig-KRSS}
\end{figure}

As an example, the Kumar-Rybko-Seidman-Stolyar (KRSS) network is
illustrated in Figure \ref{fig-KRSS}. This network, widely known as
Kumar-Seidman network and Rybko-Stolyar network in queueing network
literatures, was first studied independently by Kumar and Seidman
(1990) and Rybko and Stolyar (1992).
The KRSS network consists of two stations and four job classes.
Among the four job classes, only class 1 and 3 have external job
arrivals, i.e., $\alpha_2 = \alpha_4 = 0$. A class 1 (class 3) job
becomes a class 2 (class 4) job after its service completion at
station 1 (station 2), while a class 2 (class 4) job leaves the
system after its service completion at station 2 (station 1). The
class 4 (class 2) jobs have higher priority than class 1 (class 3)
jobs at station 1 (station 2). In particular for this network, the
parameters $C$, $P$ and $\pi$ are specified as
\begin{eqnarray}
&&C=\left(%
    \begin{array}{cccc} 1&0&0&1 \\0&1&1&0 \end{array}
    \right),~%
  P=\left(%
    \begin{array}{cccc} 0&1&0&0 \\0&0&0&0\\0&0&0&1 \\0&0&0&0 \end{array}
    \right),~\mbox{ and }%
\pi = (4,1,2,3).    \nonumber
\end{eqnarray}
With a little thought, it is direct to see that the traffic
intensity is simply
    $$\rho = (\rho_1,\rho_2)^T
    = (\alpha_1 m_1+\alpha_3 m_4,\alpha_1 m_2+ \alpha_3 m_3)^T .$$
It is well known (e.g., Chen and Zhang (2000)) that the KRSS network
is stable if and only if
    $$\rho<e \mbox{~~and~~} \alpha_1 m_2+ \alpha_3 m_4<1 . $$

The counterexample that presents a paradox in the admission
control of open multiclass queueing network is a variation of the
KRSS network. It is illustrated in Figure \ref{fig-KRSS-modi}, and
called {\it modified} KRSS network in the following.
\begin{figure}[htb]
\setlength{\unitlength}{1cm}
\begin{center}
\begin{picture}(5,4)(1,0.5)
\multiput(1,1.5)(4,0){2}{\framebox(1,2)} 
\put(3,2.55){\framebox(1,1.95)} %
\put(3,0.5){\framebox(1,1.95)} %
\multiput(0,3)(2,0){4}{\vector(1,0){1}} 
   \multiput(1.05,3)(0.2,0){5}{\line(1,0){0.1}}
   \multiput(3.05,3)(0.2,0){5}{\line(1,0){0.1}}
   \multiput(5.05,3)(0.2,0){5}{\line(1,0){0.1}}
\multiput(2,4)(2,0){2}{\vector(1,0){1}} 
   \multiput(3.05,4)(0.2,0){5}{\line(1,0){0.1}}
\multiput(1,2)(2,0){4}{\vector(-1,0){1}} 
   \multiput(1.05,2)(0.2,0){5}{\line(1,0){0.1}}
   \multiput(3.05,2)(0.2,0){5}{\line(1,0){0.1}}
   \multiput(5.05,2)(0.2,0){5}{\line(1,0){0.1}}
\multiput(5,1)(-2,0){2}{\vector(-1,0){1}} 
   \multiput(3.05,1)(0.2,0){5}{\line(1,0){0.1}}
\put(1.1,1.2){St. 1} \put(5.1,1.2){St. 2} \put(3.1,4.6){St. 3} \put(3.1,0.2){St. 4}%
\put(1.2,3.1){$m_1$} \put(5.2,3.1){$m_2$} \put(1.2,1.75){$m_4$}
\put(5.2,1.75){$m_3$} \put(3.2,2.7){$m_5$} \put(3.2,4.1){$m_7$}
\put(3.2,2.1){$m_6$} \put(3.2,0.7){$m_8$} \put(1.3,2.3){\large
$\Lambda$}\put(3.3,1.3){\large $\Lambda$} \put(5.3,2.3){\large V}
\put(0.1,3.1){$\alpha_1$}\put(6.7,1.7){$\alpha_3$}
\put(2.1,4.1){$\alpha_7$}\put(4.7,0.7){$\alpha_8$}
\put(3.3,3.3){\large V}
\end{picture}
\end{center}
\caption{Modified KRSS network} \label{fig-KRSS-modi}
\end{figure}
Compared to the original KRSS network, there are two additional
stations, namely the station 3 and 4, and four job classes, namely
the class 5, 6, 7, and 8. The class 7 (class 8) has higher priority
than the class 5 (class 6) at station 3 (station 4). The details of
the specific network parameters (${\cal J}, {\cal K}, \alpha, m, C,
P, \pi$) for this network should be obvious from the figure. For the
modified KRSS network, we have the following result.
\begin{thm} \label{thm-KRSS-modi}
Suppose $\rho<e$ and
\begin{equation}
\alpha_1 m_2 + \alpha_3 m_4 >1 \label{eq-virtual-st}
\end{equation}
in the modified KRSS network.

\noindent (1) If $m_5/(1-\alpha_7 m_7) > m_2$ and $m_6/(1-\alpha_8
m_8) > m_4$, then the queue length process $Q(\cdot)$ is positive
recurrent, and thus the expected stationary total queue length
$\bar Q < \infty$.

\noindent (2) If $m_5/(1-\alpha_7 m_7) < m_1$ and $m_6/(1-\alpha_8
m_8) < m_3$, then the queue length process $Q(\cdot)$ is
transient, and thus the expected stationary total queue length
$\bar Q = \infty$.
\end{thm}
This theorem presents a phenomenon in which reducing the arrival
rates of some job classes leads to worse performance of the
queueing network. To see this, fix all the parameters of the
modified KRSS network except $\alpha_7$ and $\alpha_8$. In
statement (1) of the above theorem, we have that
\begin{equation}
 \alpha_7 > (1-m_5/m_2)/m_7 \mbox{ and } \alpha_8 > (1-m_6/m_4)/m_8
 \label{eq-sn-110}
\end{equation}
and that the expected stationary total queue length $\bar Q$ is
finite. However, when we reduce $\alpha_7$ and $\alpha_8$ to the
case such that
\begin{equation}
 \alpha_7 < (1-m_5/m_1)/m_7  \mbox{ and } \alpha_8 < (1-m_6/m_3)/m_8,
 \label{eq-sn-120}
\end{equation}
the queue length process $Q(t)$ becomes transient and thus $\bar Q$
becomes infinite. Virtually, we will see in next section that
$\sum_{k\in {\cal K}} Q_k(t) \rightarrow \infty$ almost surely.

To gain better intuition of the paradoxical phenomenon, examine the
dynamics of the original KRSS network with no initial job (note that
the initial condition has no impact on the long term network
behavior). When a class 4 job is being served, class 1 jobs can not
move to class 2 for further service, and vice versa. From this
observation, it is not difficult to infer that classes 2 and 4 will
never be served at the same time and in effect form a virtual
station (Dai and Vande Vate 1996). Therefore, the total nominal
traffic intensity for these two classes together, i.e., the virtual
station, should not exceed one for the network to be stable. The
similar argument also yields that the KRSS network is unstable when
the nominal traffic intensity for the virtual station exceed one,
i.e., the condition (\ref{eq-virtual-st}) holds. Now consider the
modified KRSS network. The additional classes 5 and 6 act as {\it
regulators} that regulate the traffics to classes 2 and 4
respectively so as to stabilize the network. (Readers may refer to
Humes (1994) on the application of regulators to stabilize queueing
networks.) When the workloads of classes 7 and 8 are light such that
the condition (\ref{eq-sn-120}) holds, much service capacity of
stations 3 and 4 is left to classes 5 and 6 respectively and hence
the classes 5 and 6 do not hold back the traffics to avoid building
up of job queues at classes 2 and 4 respectively (cf. the case (2)
of Theorem \ref{thm-KRSS-modi}). Thus, the virtual station effect
prevails and the network is still unstable under the condition
(\ref{eq-virtual-st}). However, when the workloads of classes 7 and
8 are heavy enough such that the condition (\ref{eq-sn-110}) holds,
the service for lower priority classes 5 and 6 is in effect slowed
down and the traffics to classes 2 and 4 are held back (cf. the case
(1) of Theorem \ref{thm-KRSS-modi}). Consequently, there would not
be large buildup of queues at classes 2 and 4, and these two classes
will not mutually block their services. Finally, the virtual station
effect is avoided and the modified KRSS network is thus stabilized.
The above argument will be made rigorous in the proof of Theorem
\ref{thm-KRSS-modi} in next section.

\begin{figure}[htbp]
 \scalebox{0.4} {\includegraphics[angle=270]{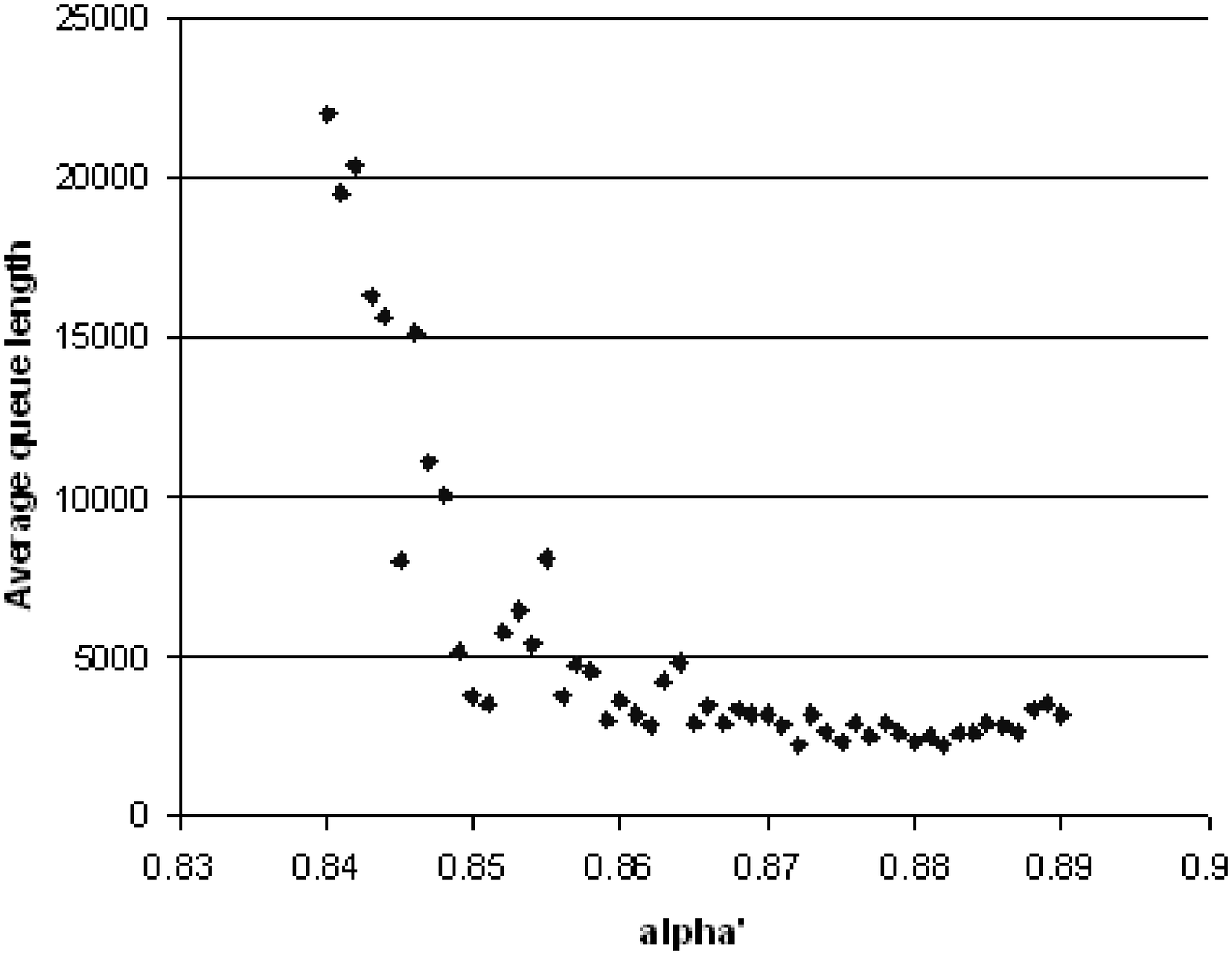}}
 \caption{Simulation for modified KRSS network}
 \label{fig-sim-KRSS}
 \end{figure}

Concerning the above paradoxical phenomenon, a subtle question is
whether this counter-intuitive phenomenon is just due to
pathological jumps in the network performance. To post this question
in more details, we take for the moment that $\alpha_1 = \alpha_3 =
1$, $\alpha_7 = \alpha_8 = \alpha'$, $m_1 = m_3 = 0.2$, $m_2 = m_4 =
0.6$, $m_5 = m_6 = 0.1$, $m_7 = m_8 = 1$.
 Then, let $\alpha'$ varies, say, from $8/9$ down to $1/3$, and thus
$m_5/(1-\alpha_7 m_7)$ and $m_6/(1-\alpha_8 m_8)$ both vary from
$0.9$ (which is greater than $m_3$ and $m_6$) to $0.15$ (which is
less than $m_3$ and $m_6$).
 Based on Theorem \ref{thm-KRSS-modi}, the expected stationary total queue
length $\bar Q$ is finite when $\alpha'$ is $8/9$, but it becomes
worse, i.e., $\bar Q=\infty$, when $\alpha'$ is reduced to $1/3$.
 Now, the subtle questions are as follows. Is this performance
degradation upon reducing arrival rate $\alpha'$ simply due to a
jump from a stable to an unstable network at a critical point of
$\alpha'$ when it varies from $8/9$ to $1/3$? Is the performance
$\bar Q$ still an increasing function of the arrival rate $\alpha'$
within any interval of $\alpha'$ where the network is stable and its
expected total queue length $\bar Q$ is finite?
 It is not obvious how to eliminate this possible
pathological situation theoretically. However, our simulation
results illustrated in Figure \ref{fig-sim-KRSS} indicate that the
average total queue length $\bar Q$ is an decreasing function of
$\alpha'$ within some intervals of $\alpha'$ (i.e., the interval
$[0.84, 0.89]$ in our simulation) where $Q(t)$ is stable. In words,
the network performance is improved {\it continuously} when more
jobs are admitted to the system within certain range of job arrival
rates.

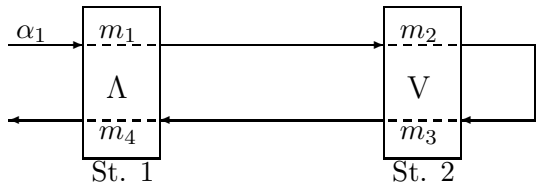
\begin{figure}[htb]
\setlength{\unitlength}{1cm}
\begin{center}
\begin{picture}(5,3)(1,1.5)
\multiput(1,1.5)(4,0){2}{\framebox(1,2)} 
\put(0,3){\vector(1,0){1}} 
\put(2,3){\vector(1,0){3}} 
\put(6,3){\line(1,0){1}} 
\put(7,3){\line(0,-1){1}} 
\multiput(1,2)(6,0){2}{\vector(-1,0){1}} 
\put(5,2){\vector(-1,0){3}} 
   \multiput(1.05,3)(0.2,0){5}{\line(1,0){0.1}}
   \multiput(5.05,3)(0.2,0){5}{\line(1,0){0.1}}
   \multiput(1.05,2)(0.2,0){5}{\line(1,0){0.1}}
   \multiput(5.05,2)(0.2,0){5}{\line(1,0){0.1}}
\put(1.1,1.2){St. 1} \put(5.1,1.2){St. 2}
\put(0.1,3.1){$\alpha_1$} \put(1.2,3.1){$m_1$}
\put(5.2,3.1){$m_2$} \put(1.2,1.75){$m_4$} \put(5.2,1.75){$m_3$}
\put(1.3,2.3){\large $\Lambda$} \put(5.3,2.3){\large V}
\end{picture}
\end{center}
\caption{LK network} \label{fig-LK}
\end{figure}

\begin{figure}[htb]
\setlength{\unitlength}{1cm}
\begin{center}
\begin{picture}(5,3.5)(1,1.5)
\multiput(1,1.5)(4,0){2}{\framebox(1,2)} 
\put(3,2.5){\framebox(1,2)} %
\multiput(0,3)(2,0){3}{\vector(1,0){1}} 
\put(6,3){\line(1,0){1}} 
\put(7,3){\line(0,-1){1}} 
   \multiput(1.05,3)(0.2,0){5}{\line(1,0){0.1}}
   \multiput(3.05,3)(0.2,0){5}{\line(1,0){0.1}}
   \multiput(5.05,3)(0.2,0){5}{\line(1,0){0.1}}
\put(5.5,4){\vector(-1,0){1.5}} 
\put(3,4){\line(-1,0){3}} 
\put(0,4){\line(0,-1){1}} 
   \multiput(3.05,4)(0.2,0){5}{\line(1,0){0.1}}
\multiput(1,2)(6,0){2}{\vector(-1,0){1}} 
\put(5,2){\vector(-1,0){3}} 
   \multiput(1.05,2)(0.2,0){5}{\line(1,0){0.1}}
   \multiput(5.05,2)(0.2,0){5}{\line(1,0){0.1}}
\put(1.1,1.2){St. 1} \put(5.1,1.2){St. 2} \put(3.1,4.6){St. 3}
\put(1.2,3.1){$m_1$} \put(5.2,3.1){$m_2$} \put(1.2,1.75){$m_4$}
\put(5.2,1.75){$m_3$} \put(3.2,2.7){$m_5$} \put(3.2,4.1){$m_6$}
\put(1.3,2.3){\large $\Lambda$} \put(5.3,2.3){\large V}
\put(4.7,4.1){$\alpha_6$} \put(3.3,3.3){\large V}
\end{picture}
\end{center}
\caption{A modified LK network} \label{fig-LK-modi}
\end{figure}
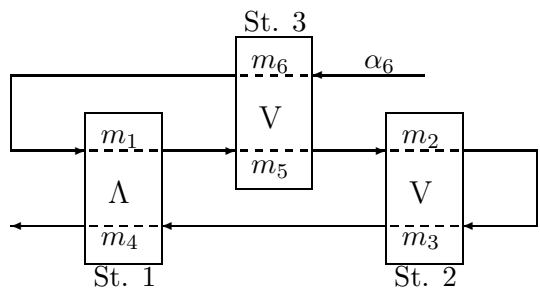

Another counterexample that gives different perspective on the
paradox in admission control is related to the Lu-Kumar (LK)
network, which was first studied by Lu and Kumar (1991) and is
illustrated in Figure \ref{fig-LK}. We omit the detailed
description of this network, which should be clear from the its
comparison with the KRSS network.
This counterexample is a variation of the LK network, called {\it
modified} LK network in this paper, and is illustrated in Figure
\ref{fig-LK-modi}.
For the modified LK network, we study some special instances (for
convenience) and summarize the counter-intuitive phenomenon in the
following theorem.
\begin{thm} \label{thm-LK-modi}
Consider the modified LK network with $m = (0.1, 0.6, 0.1, 0.6, 0.7,
0.027)^T$.

\noindent (1) If $\alpha_6 = 1.37$, then the queue length process
$Q$ is positive recurrent, and thus $\bar Q < \infty$.

\noindent (2) If $\alpha_6 = 1$, then the queue length process $Q$
is transient, and thus $\bar Q = \infty$.
\end{thm}
This theorem presents a situation in which, when the arrival rate
$\alpha_6$ drops from $1.37$ to $1$, the performance becomes worse.
Similar to the simulation for the modified KRSS network, our
simulation result also supports that  for the modified LK network
the average total queue length $\bar Q$ would be a decreasing
function of $\alpha_6$ within some intervals of $\alpha_6$ where
$Q(t)$ is stable.
 In contrast to the modified KRSS network, a
special feature of the modified LK network is that there is only one
external arrival and this arrival is controllable. On the other
hand, if we fix the rate $\alpha_6$ of the unique external arrival
and vary the service times $m_k$, $k =1,\cdots,6$, in proportion,
then we recover an example for the paradox on service control. That
is, increasing service capacity may also worsen the system
performance, since reducing the service times $m_k$ in proportion
(i.e., increasing the service capacity) is equivalent to reducing
the external arrival $\alpha_6$ in the modified LK network by
changing the time scale suitably.

\section{Multiclass Fluid Network Model
and Proof of Theorem \ref{thm-KRSS-modi}}

In this section, we provide the proof of Theorem
\ref{thm-KRSS-modi}, while the proof of Theorem \ref{thm-LK-modi} is
omitted since it is similar to the former one. We employ the fluid
model approach in the proof.
The development of this approach was inspired by the studies of some
counterexamples in Kumar and Seidman (1990), Rybko and Stolyar
(1992) and Bramson (1994), etc., where the multiclass queueing
networks are not stable even when the traffic intensity of each
station in the network is less than one. An elegant result of the
fluid model approach was proposed first in Rybko and Stolyar (1992)
and then generalized and refined by Dai (1995), Chen (1995), Dai and
Meyn (1995), Stolyar (1995) and Bramson (1998). It states that a
queueing network is stable if its corresponding fluid network model
is stable. Partial converse to this result is also given in Meyn
(1995), Dai (1996) and Puhalskii and Rybko (2000). To quote these
results to prove our theorems, we first present in next subsection a
multiclass fluid network model corresponding to the multiclass
queueing network model described in Section \ref{sec-main}.

\subsection{A Multiclass Fluid Network Model}

Parallel to the queueing network (${\cal J}, {\cal K}, \alpha, m, C,
P, \pi$), a corresponding fluid network model is obtained
intuitively by replacing the discrete jobs in the queueing network
with continuous fluids. Specifically, the fluid network also
consists of $J$ stations (buffers) indexed by  $j \in {\cal J} =
\{1,...,J \}$, serving $K$ fluid (job) classes indexed by $k \in
{\cal K}= \{ 1,...,K\}$. A fluid class is served exclusively at one
station, but one station may serve more than one fluid classes. As
in the queueing network, $\sigma(\cdot)$ denotes a many-to-one
mapping from $\cal K$ to $\cal J$, with $\sigma(k)$ indicating the
station at which a class $k$ fluid is served.  A class $k$ fluid may
flow exogenously into the network at rate $\alpha_k(\ge0)$, then it
is served at station $\sigma(k)$, and after being served, a fraction
$p_{k\ell}$ of fluid turns into a class $\ell$ fluid, $\ell\in{\cal
K}$, and the remaining fraction, $1-\sum_{\ell=1}^Kp_{k\ell}$ flows
out of the network.   When station $\sigma(k)$ devotes its full
capacity to serving class $k$ fluid (assuming that it is available
to be served), it generates an outflow of class $k$ fluid at rate
$\mu_k(>0)$, $k\in{\cal K}$.  Let $\alpha=(\alpha_k)$ and call it
the exogenous {\it inflow (arrival) rate} (vector). Let
$\mu=(\mu_k)$ and call it the {\it service rate} (vector).  We call
$K\times K$ substochastic matrix $P=(p_{k\ell})$ the {\it flow
transition matrix}. Corresponding to the open queueing network
described in the last section, we consider an open fluid network.
That is, we also assume that  matrix $P$ satisfies the condition
(\ref{eq-sn-80}).
Among classes, fluid follows a priority service discipline, which is
again described by a one-to-one mapping $\pi$ from $\{1,...,K\}$
onto itself. Specifically, a class $k$ has priority over a class
$\ell$ if $\pi(k)<\pi(\ell)$. We adopt the following notation from
the description of the multiclass queueing network model, i.e., $C$,
$\lambda$, $\beta$, $\rho$, $M$, and $D$.

To describe the dynamics of the fluid network, we introduce the
$K$-dimensional fluid level process ${\bar Q}=\{{\bar
Q}(t),t\ge0\}$, whose $k$th component ${\bar Q}_k(t)$ denotes the
fluid level of class k at time $t$; the $K$-dimensional time
allocation process ${\bar T}=\{{\bar T}(t), t\ge0\}$, whose $k$th
component ${\bar T}_k(t)$ denotes the total amount of time that
station $\sigma (k)$ has devoted to serving class $k$ fluid during
the time interval $[0,t]$; and the $K$-dimensional unused capacity
process ${\bar Y}=\{{\bar Y}(t),t\ge0\}$, whose $k$th component
${\bar Y}_k(t)$ denotes the  (cumulative) unused capacity of station
$\sigma (k)$ during the time interval $[0,t]$ after serving all
classes at station $\sigma (k)$ which have a priority no less than
class $k$ (including class $k$). Let
    $$H_k = \{ \ell: \sigma(\ell)=\sigma(k),
    \pi(\ell) \leq \pi(k) \}$$
be the set of indices for all classes that are served at the same
station as class $k$ and have a priority no less than that of class
$k$. Note that $k\in H_k$ by definition. Then the dynamics of the
fluid network model can be described as follows.
\begin{eqnarray}
&& \bar Q(t)=\bar Q(0) + \alpha t - (I-P') D \bar T(t)\ge0, \label{eq-fn-22-new}\\
&& \bar T(\cdot)\mbox{ is nondecreasing with } \bar T(0)=0, \label{eq-fn-23-new}\\
&& \bar Y_k(t) = t-\sum_{\ell\in H_k} \bar T_\ell(t) \mbox{ is
nondecreasing},\quad k\in{\cal K},\label{eq-fn-22}\\
&& \int ^{\infty}_0 \bar Q_k(t)d \bar Y_k (t) =0, \qquad k\in{\cal
K}. \label{eq-fn-26}
\end{eqnarray}
The relation (\ref{eq-fn-22-new}) is the flow balance relation;
its $k$th coordinate reads as,
    $$Q_k(t)=Q_k(0) + \alpha_k t + \sum_{\ell=1}^K p_{\ell k} \mu_\ell
    T_\ell(t) - \mu_k T_k(t) \ge0,~~k=1,\cdots,K .$$
The equation (\ref{eq-fn-22-new}) is nothing but the equivalent
relation between the time allocation process $T(\cdot)$ and the
unused capacity process $Y(\cdot)$. The relation (\ref{eq-fn-26})
specifies both the work-conserving condition and the priority
discipline; in words, for each $k$, the relation (\ref{eq-fn-26})
means that at any time $t$,  there could be some positive remaining
capacity (rate) for serving those classes at station $\sigma(k)$
having a strictly lower priority than class $k$, only when the fluid
levels of all classes in $H_k$ (having a priority no less than $k$)
are zero. Particularly, for each lowest fluid class $k$ at station
$j=\sigma(k)$, the relation (\ref{eq-fn-26}) specifies the
work-conserving condition for station $j$, which implies that
station $j$ can not be idle if the total fluid level
($\sum_{\ell:\sigma(\ell)=j} \bar Q_\ell(t)$) in station $j$ is
positive at any time $t\ge 0$.

We shall refer to this network as fluid network (${\cal J}, {\cal
K}, \alpha, m, C, P, \pi$). For the fluid network (${\cal J}, {\cal
K}, \alpha, m, C, P, \pi$), A pair $(\bar Q,\bar T)$ (or
equivalently ($\bar Q, \bar Y$)) is said to be a {\it fluid
solution} if they jointly satisfy
(\ref{eq-fn-22-new})-(\ref{eq-fn-26}). For convenience, we also call
$\bar Q$ a fluid solution if there is a $\bar T$ such that the pair
($\bar Q, \bar T$) is a fluid solution. The fluid network (${\cal
J}, {\cal K}, \alpha, m, C, P, \pi$) is said to be {\it stable} if
there is a time $\tau\ge0$ such that $\bar Q(\tau+\cdot)\equiv0$ for
any fluid solution $\bar Q$ with $||\bar Q(0)||=1$; and it is said
to {\it weakly stable} if $\bar Q(\cdot)\equiv0$ for any fluid
solution $\bar Q$ with $\bar Q(0)=0$. A well-known property we will
use later in this paper is that the processes $\bar Q$, $\bar Y$,
and $\bar T$ are Lipschitz continuous, and hence are differentiable
almost everywhere on $[0 , \infty )$. We summarize some known
stability results on the relation between the queueing network model
and its corresponding fluid network model, which are used in the
proof of Theorem \ref{thm-KRSS-modi}.
\begin{thm} \label{thm-QN-FN}
Consider the queueing network (${\cal J}, {\cal K}, \alpha, m, C,
P, \pi$).

\noindent (1) If the corresponding fluid network (${\cal J}, {\cal
K}, \alpha, m, C, P, \pi$) is stable, then the queue length
process $Q$ is positive recurrent.

\noindent (2) If the corresponding fluid network (${\cal J}, {\cal
K}, \alpha, m, C, P, \pi$) is not weakly stable, then the queue
length process $Q$ is transient.
\end{thm}
Readers are referred to Chen and Yao (2001) and Dai (1996) for
elementary proofs of the two conclusions respectively.

\subsection{Proof of Theorem}

\noindent {\bf Proof of Theorem \ref{thm-KRSS-modi} (1):} According
to Theorem \ref{thm-QN-FN} (1), it is sufficient to show that the
fluid network model corresponding to the modified KRSS queueing
network, called the modified KRSS fluid network below, is stable. As
an instance of the fluid network model described in
(\ref{eq-fn-22-new})-(\ref{eq-fn-26}), the dynamics of the modified
KRSS fluid network can be detailed as follows.
\begin{eqnarray}
&& \bar Q_k(t)=\bar Q_k(0) + \alpha_k t - \mu_k \bar T_k(t)\ge0,
~~~~k=1,3,7,8, \label{fn-modi-KRSS-1}\\
&& \bar Q_k(t)=\bar Q_k(0) + \mu_\ell \bar T_\ell(t) - \mu_k \bar
T_k(t)\ge0,
~~~~(k,\ell)=(5,1),(2,5), (6,3) , (4,6), \label{fn-modi-KRSS-1b}\\
&& \bar T_k(\cdot)\mbox{ is nondecreasing with } \bar T_k(0)=0,
~~~~k=1,\cdots,8, \label{fn-modi-KRSS-2}\\
&& \bar Y_k(t) = t- \bar T_k(t) \mbox{ is
nondecreasing},\quad k=4,2,7,8,\label{fn-modi-KRSS-3}\\
&& \bar Y_k(t) = t- \bar T_\ell(t) - \bar T_k(t) \mbox{ is
nondecreasing},\quad (k,\ell)=(1,4),(3,2),(5,7),(6,8),\label{fn-modi-KRSS-3b}\\
&& \int ^{\infty}_0 \bar Q_k(t)d \bar Y_k (t) =0, \qquad
k=1,\cdots,8. \label{fn-modi-KRSS-4}
\end{eqnarray}

We prove the stability the modified KRSS fluid network in three
steps. First, we prove that there exists a time $\tau_1 \ge 0$
such that
\begin{eqnarray}
&& \bar Q_7(t)=\bar Q_8(t)=0 \mbox{ for any } t\ge \tau_1.
\label{proof-KRSS-1}
\end{eqnarray}
If $\dot {\bar Q}_7(t) >0$, then we have
\begin{eqnarray}
&& \dot{\bar Y}_7(t) =0 \label{proof-KRSS-5}
\end{eqnarray}
by condition (\ref{fn-modi-KRSS-4});
\begin{eqnarray}
&& \dot{\bar T}_7(t) =1 \label{proof-KRSS-10}
\end{eqnarray}
by equations (\ref{fn-modi-KRSS-3}) and (\ref{proof-KRSS-5}); and at
last
\begin{eqnarray}
&& \dot{\bar Q}_7(t) = \alpha_7 - \mu_7 \label{proof-KRSS-15}
\end{eqnarray}
by equations (\ref{fn-modi-KRSS-1}) and (\ref{proof-KRSS-10}). Note
that the condition $\rho < e$ implies $\alpha_7 < \mu_7 $. Let
$\tau_1' = { \dot{\bar Q}_7(0)} / {(\mu_7 - \alpha_7)}$. Then, we
have
\begin{eqnarray}
&& \bar Q_7(t)=0 \mbox{ for any } t\ge \tau_1'.
\label{proof-KRSS-20}
\end{eqnarray}
Using the similar argument, we have
\begin{eqnarray}
&& \bar Q_7(t)=0 \mbox{ for any } t\ge \tau_1'' = \frac{ \dot{\bar
Q}_7(0)}{\mu_8 - \alpha_8}. \label{proof-KRSS-25}
\end{eqnarray}
Letting $\tau_1 = \max({1}/{(\mu_7 - \alpha_7)}, {1}/{(\mu_8 -
\alpha_8)} )$, we have that $\tau_1 \ge \max(\tau_1',\tau_1'')$
under the assumption $||\bar Q(0)|| = 1$. Now, the conclusions
(\ref{proof-KRSS-20}) and (\ref{proof-KRSS-25}) lead to the claim
(\ref{proof-KRSS-1}).

Next, we prove that there exists a time $\tau_2 \ge \tau_1 $ such
that
\begin{eqnarray}
&& \bar Q_4(t)=\bar Q_2(t)=0 \mbox{ for any } t\ge \tau_2.
\label{proof-KRSS-30}
\end{eqnarray}
Under the condition (\ref{proof-KRSS-1}), we have $\dot{\bar
Q}_7(t)=\dot{\bar Q}_8(t)=0$, and then $\dot{\bar T}_7(t)=\alpha_7
m_7$ and $\dot{\bar T}_8(t)=\alpha_8 m_8$ for all time $t\ge
\tau_1$. Combined with (\ref{fn-modi-KRSS-3b}), this gives rise to
\begin{eqnarray}
&& \dot{\bar Y_6}(t)= 1- \dot{\bar T_6}(t)-\dot{\bar T_8}(t)\ge 0,
\mbox{~~and } \nonumber \\
&& \dot{\bar T_6}(t) \le 1 - \dot{\bar T_8}(t) = 1 - \alpha_8 m_8,
 \mbox{~~~~for any } t\ge \tau_1. \label{proof-KRSS-35}
\end{eqnarray}
Then, we have
\begin{eqnarray}
&& \dot{\bar Q_4}(t)= \mu_6 \dot{\bar T_6}(t)- \mu_4 \dot{\bar
T_4}(t) \le \mu_6(1-\alpha_8 m_8) - \mu_4 <0  \nonumber
\end{eqnarray}
for any $t\ge\tau$, where the last inequality is implied by the
assumption that $m_6/(1-\alpha_8 m_8) > m_4$. Let $\tau_2' =
 \dot{\bar Q}_4(\tau_1)/ ( \mu_4 - \mu_6(1-\alpha_8 m_8) )$.
Then, we have
\begin{eqnarray}
&& \bar Q_4(t)=0 \mbox{ for any } t\ge \tau_2'.
\label{proof-KRSS-40}
\end{eqnarray}
Similarly, we have
\begin{eqnarray}
&& \bar Q_2(t)=0 \mbox{ for any } t\ge \tau_2'' = \frac{ \dot{\bar
Q}_2(\tau_1)}{\mu_2 - \mu_5(1-\alpha_7 m_7)}.
\label{proof-KRSS-45}
\end{eqnarray}
Let
    $$\tau_2 = \max(\frac{1 + \Delta \tau_1}{\mu_4 -
    \mu_6(1-\alpha_8 m_8)},\frac{1+ \Delta \tau_1}{\mu_2 -
    \mu_5(1-\alpha_7 m_7)} )$$
with $\Delta$ being the Lipschitz constant for the fluid level
process $\bar Q(t)$. Then we have that $\tau_2 \ge
\max(\tau_2',\tau_2'')$, noting that $||\bar Q(\tau_1)|| \le ||\bar
Q(\tau_1)|| + M \tau_1 \le 1 + M \tau_1$. Now, the conclusions
(\ref{proof-KRSS-40}) and (\ref{proof-KRSS-45}) imply the claim
(\ref{proof-KRSS-30}).

Finally, we prove that there exists a time $\tau \ge \tau_2
(\ge0)$ such that
\begin{eqnarray}
&& \bar Q_k(t)=0 \mbox{~~ for } k=1,3,5,6 \mbox{ and } t\ge \tau,
\label{proof-KRSS-50}
\end{eqnarray}
which together with equations (\ref{proof-KRSS-1}) and
(\ref{proof-KRSS-30}) implies
\begin{eqnarray}
&& \bar Q(t)=0 \mbox{~~ for } t\ge \tau. \nonumber
\end{eqnarray}
Let
\begin{eqnarray}
\bar W_1(t) &:=& m_1 \bar Q_1(t) +  m_4 (\bar Q_3(t) + \bar
Q_6(t)) \nonumber \\
&=& (\alpha_1 m_1 + \alpha_3 m_4)t - (\bar T_1(t) + \bar T_4(t)) , \nonumber \\
\bar W_2(t) &:=& m_3 \bar Q_3(t) +  m_2 (\bar Q_1(t) + \bar
Q_5(t)) \nonumber \\
&=& (\alpha_1 m_2 + \alpha_3 m_3)t - (\bar T_2(t) + \bar T_3(t)), \nonumber \\
\bar W_3(t) &:=& m_5 (\bar Q_1(t) + \bar Q_5(t))
= \alpha_1 m_5 t - \bar T_5(t) , \nonumber \\
\bar W_4(t) &:=& m_6 (\bar Q_3(t) + \bar Q_6(t)) = \alpha_3 m_6 t
- \bar T_6(t) , \nonumber
\end{eqnarray}
for $t\ge \tau_2$. Here $\bar W_i(t)$ ($i=1,2,3,4,$) can be
explained as the immediately workload for station $i$ implied in
the system at time $t$. Define
\begin{eqnarray}
&& f_1(t) := m_6 \bar W_1(t), ~f_2(t) = m_5 \bar W_2(t),\nonumber\\
&& f_3(t) := m_2 \bar W_3(t),~ f_4(t) = m_4 \bar W_4(t).\nonumber
\end{eqnarray}
Then, it is direct to verify that, for $t\ge \tau_2$,
\begin{eqnarray}
&& \dot f_i(t) <0 \mbox{ if } \bar Q_i(t)>0 , \mbox{ for }
i=1,2,3,4, \nonumber
\end{eqnarray}
and
\begin{eqnarray}
&& f_1(t) \le f_4(t) \mbox{ if } \bar Q_1(t)=0 ,\nonumber\\
&& f_2(t) \le f_3(t) \mbox{ if } \bar Q_3(t)=0 ,\nonumber\\
&& f_3(t) \le f_2(t) \mbox{ if } \bar Q_5(t)=0 ,\nonumber\\
&& f_4(t) \le f_4(t) \mbox{ if } \bar Q_6(t)=0 .\nonumber
\end{eqnarray}
Now applying the piecewise linear Lyapunov function approach for the
multiclass fluid network model described in Theorem 3.1 of Chen and
Ye~(2002), we obtain the conclusion (\ref{proof-KRSS-50}).
 \hfill Q.E.D.\medskip 

\noindent {\bf Proof of Theorem \ref{thm-KRSS-modi} (2):} ~~%
According to Theorem \ref{thm-QN-FN} (2), we need to show that the
modified KRSS fluid network is not weakly stable. Similar to the
above proof of the claim (1), it is not difficult to show that there
exists a time $\tau_1 \ge 0$ such that
\begin{eqnarray}
&& \bar Q_7(t)=\bar Q_8(t)=0 \mbox{ for any } t\ge \tau_1
\nonumber
\end{eqnarray}
since classes $7$ and $8$ fluids have higher priorities at stations
$7$ and $8$ respectively; and then that there exists a time $\tau_2
\ge \tau_1 $ such that
\begin{eqnarray}
&& \bar Q_5(t)=\bar Q_6(t)=0 \mbox{ for any } t\ge \tau_2 \nonumber
\end{eqnarray}
since the remaining service capacity for classes $5$ and $6$ fluids
is greater than that for class $1$ and $3$ fluids. Thus, the
modified KRSS fluid network is reduced to the well known KRSS fluid
network, which is not weakly stable under the condition
(\ref{eq-virtual-st}). \hfill Q.E.D. 


\section{Discussion and Concluding Remark}

We have presented a paradox for the admission control for the
multiclass queueing network with differentiated service in this
paper. This paradox is, to our knowledge, the first one of the kind,
which is complement to the existing ones on the service rate control
and the routing control.

The models, as well as the admission control and the differentiated
service, studied in the paper are simplified and idealized models of
practical systems. Take the semiconductor production as the example.
The production line may consist tens of processing stations
(machines), and parts may require tens or even hundreds of stages of
processing by the stations. The admission control may model the {\it
central control} on whether to accept the external order, while the
differentiated priority for jobs at each station could be due to the
{\it local control} on scheduling jobs. In addition, machines may be
subject to random failures and need set-up time when changing from
processing a class of jobs to another.
 it would not be surprising that the paradoxical phenomenon in the
admission control exists in such a complex and practical system, as
it exists even in the simplified and idealized network models
presented in this paper.
 Therefore, the detection of and the remedy to such a
paradoxical phenomenon would be interesting future research topics.


\end{document}